# SLE($\kappa, \rho$) MARTINGALES AND DUALITY

By Julien Dubédat

*Université Paris-Sud*

Various features of the two-parameter family of Schramm–Loewner evolutions SLE($\kappa, \rho$) are studied. In particular, we derive certain restriction properties that lead to a "strong duality" conjecture, which is an identity in law between the outer boundary of a variant of the SLE($\kappa$) process for $\kappa \geq 4$ and a variant of the SLE($16/\kappa$) process.

**1. Introduction.** Stochastic (or Schramm) Loewner evolutions (SLEs) are stochastic increasing families of plane compact sets. Loewner showed how to parametrize an increasing family of compact sets ("hulls") in a plane domain with a real-valued continuous function (under a "local growth" condition). In other words, Loewner equations transform a real path into an increasing family of hulls. SLEs, which were first introduced by Schramm [21], are basically the image of the Wiener measure under this transformation. It turns out that the measures on compact sets obtained in this way have very different properties according to the speed $\kappa$ of the driving Brownian motion.

Since the SLEs give probability laws on hulls that have built-in conformal invariance properties, they are the only possible candidates for the scaling limits of various critical plane discrete models, which are conformally invariant in the scaling limit. The cases for which the convergence to the scaling limit has been proved are uniform spanning trees (UST), loop-erased random walks (LERW) and critical percolation. Smirnov [23] proved that the scaling limit of critical percolation clusters on the triangular lattice is described by SLE(6). Lawler, Schramm and Werner [14], among numerous results on SLE, proved that the scaling limit of the UST Peano curve (resp. the LERW) is SLE(8) [resp. SLE(2)]. For other critical models, conformal invariance is conjectured but not proved: Double domino tiling paths are believed to converge to SLE(4) (see [20]), critical FK percolation cluster interfaces are conjectured to converge to SLE($\kappa$) (the $q$ parameter of the FK percolation and the $\kappa$ parameter of SLE are linked by the conjectural relation $-\sqrt{q}/2 = \cos(4\pi/\kappa)$, see [20]) and there is some evidence that self-avoiding walks should converge to SLE(8/3) [15].

There are two main variants of SLEs: chordal SLEs, which depend on a domain and two points on the boundary (or prime ends), and radial SLEs, which depend on a domain, an inner point and a point on the boundary. We are mainly interested in chordal SLE, but this is not so restrictive, since the radial and the chordal constructions are "equivalent" in some appropriate sense for small enough times; see [12], Proposition 4.2.

As we already mentioned, the value of the $\kappa$ parameter has a big influence on the geometric properties of the SLE. When $\kappa \leq 4$, SLEs are a.s. simple paths [20]. If $\kappa > 4$, this is no longer the case, but SLEs are generated by a continuous path, the trace, that can have double points, but cannot cross its past (see [20]). When $\kappa \geq 8$, the trace becomes space-filling. The "phase transition" at $\kappa = 4$ separates SLEs that are simple paths from SLEs that have a nontrivial boundary (for finite times, since SLEs eventually swallow the whole space).

Conjectures on the Hausdorff dimensions of these SLE paths and their outer boundaries prompted Duplantier and others to formulate the following concept.

CONJECTURE 1 (Duality for SLE). *When $\kappa > 4$, the boundary of a SLE($\kappa$) looks locally like a SLE($16/\kappa$).*

We record a very loose formulation on purpose, since actually getting an accurate statement is not straightforward. Note that it is not very difficult to guess the dimension of the SLE paths and their outer boundary by roughly evaluating the probability that a given point is on the $\varepsilon$ neighborhood of these sets. Assuming that this guess is correct and that the outer boundary of a SLE (for $\kappa \geq 4$) also looks like a SLE, it is then natural to conjecture that it should look like a SLE($16/\kappa$) which is the only one with the appropriate Hausdorff dimension.

In fact, the dimension of the Hausdorff dimension of a SLE($\kappa'$) was proved [2] to be $1 + \kappa'/8$, while the result for the Hausdorff dimension of the outer boundary of a SLE($\kappa$) for $\kappa > 4$ is still conjectural. Hence, a proof of the duality conjecture would, in particular, imply that when $\kappa > 4$, the dimension of the outer boundary of a SLE($\kappa$) has dimension $1 + 2/\kappa$. A direct proof of this fact might exist.

This duality conjecture actually was proved in the two special cases $\kappa = 8$ and $\kappa = 6$. For $\kappa = 8$, $16/\kappa = 2$, the result follows by, respectively, identifying the SLE(8) with the scaling limit of the UST Peano curve and the SLE(2) with the scaling limit of the LERW. An exact relationship (Pemantle [18] or Wilson's [27] algorithm) is known between these two discrete models, which leads to a relationship in the scaling limit, using the convergence derived in [14].

The relationship for $\kappa = 6$ was established in [16] in a way that is closely related to the approach of our paper. In [16], all the random subsets of a domain that satisfy the "conformal restriction property" are described (we briefly recall this property in Section 3; loosely speaking, the laws of these random sets are invariant under a certain semigroup of conformal transformations). It turns out that the outer boundaries of these sets are all variants of the SLE(8/3) process: the SLE(8/3, $\rho$) processes that all "look locally like" the SLE(8/3) process. They can be viewed as a SLE(8/3) process with an additional drift away from (or toward) one part of the boundary.

Other sets can directly be shown to satisfy this restriction property: conditioned Brownian motions and conditioned SLE(6) processes. This yields a description of their outer boundary in terms of a SLE(8/3, $\rho$) process [16]. The conditioned SLE(6) can be understood as follows: It is a SLE(6) from $a$ to $b$ in a domain $D$ that is conditioned "not to hit one of the two arcs between $a$ and $b$." Equivalently, consider critical percolation in a domain $D$ and condition it in such a way that no cluster touches both boundary arcs between the two boundary points $a$ and $b$. Then the conditioned SLE(6) is the scaling limit of the exploration process (exploring the boundary of the clusters attached to one part of the boundary). It turns out (see [16]) that this conditioned process is a SLE(6, 2) process.

Another way to construct the random sets that satisfy the restriction property that was pointed out in [16] is to start with a SLE($\kappa$) process for $\kappa < 8/3$ and add to this process a certain density of Brownian loops. Further properties of these Brownian loops (and the Brownian loop soup) were studied in [17]. This construction is also related to representation theory, as pointed out in [8].

In the present paper, largely based on ideas in [16], we investigate some natural generalizations of certain "restriction formulas" introduced there. In particular, we see that the properties that were derived for SLE($\kappa, 0$) (the fact that adding Brownian loops generates a set that satisfies the restriction property) on the one hand, and for SLE(8/3, $\rho$) (their relationship with the restriction measures) can be generalized to SLE($\kappa, \rho$) processes. In particular, adding a certain loop soup to SLE($\kappa, \rho$) processes for $\kappa \leq 8/3$ gives yet other ways to construct the random sets that satisfy the conformal restriction property (this was derived in the case $\kappa = 8/3$ or $\rho = 0$ in [16]). Moreover, the same computation shows that when $\kappa \geq 4$, the process SLE($\kappa, \kappa - 4$) has some special features. In particular, we prove an identity in law (when $\kappa \geq 4$ and $\kappa' = 16/\kappa \leq 4$) between the following hulls:

- The hull obtained when is a certain loop soup is added to SLE($\kappa, \kappa - 4$).

- The hull obtained when the same loop soup is added to the image of the $\mathrm{SLE}(\kappa', (\kappa'-4)/2)$ process under symmetry with respect to the imaginary axis.

This leads to the "global duality" conjecture that the outer boundary of the symmetric image of $\mathrm{SLE}(\kappa, \kappa - 4)$ is the $\mathrm{SLE}(\kappa', (\kappa'-4)/2)$ curve. We also prove that $\mathrm{SLE}(\kappa, \kappa-4)$ can be also viewed as a $\mathrm{SLE}(\kappa)$ process conditioned not to intersect one part of the boundary (which was also the case when $\kappa = 6$).

The next step to understand is a path decomposition that gives the conditional law of a $\mathrm{SLE}(\kappa)$ for $\kappa > 4$, given its outer boundary. We investigate some aspects of this question, based on restriction formulas, that lead to a stronger duality conjecture on the (local) law of the law of the hull (of the SLE), given the boundary.

This paper is organized as follows: The next two sections are an overview of $\mathrm{SLE}(\kappa)$ and $\mathrm{SLE}(\kappa, \rho)$ processes, and of the restriction formalism introduced in [16]. Section 4 presents natural generalizations of some of these results when $\kappa \neq 8/3$. In Section 5, we begin the study of the remarkable $\mathrm{SLE}(\kappa, \kappa - 4)$ processes, $\kappa > 4$, including a path decomposition. To make this decomposition more explicit, we are led to define generalized $\mathrm{SLE}(\kappa, \underline{\rho})$ processes in Section 6, before finally giving the stronger duality conjecture.

**2. Chordal SLE and $\mathrm{SLE}(\kappa, \rho)$ processes.** We first briefly recall the definition of chordal SLE in the upper half-plane $\mathbb{H}$ going from 0 to $\infty$ (see, e.g., [11, 20, 25] for more details). For any $z \in \mathbb{H}$, $t \geq 0$, define $g_t(z)$ by $g_0(z) = z$ and

$$\partial_t g_t(z) = \frac{2}{g_t(z) - W_t},$$

where $(W_t)$ is a continuous real-valued process. This ordinary differential equation (ODE) is well defined up to a random time $\tau_z$. Define the hull $K_t$ as

$$K_t = \overline{\{z \in \mathbb{H} : \tau_z < t\}}.$$

The family $(K_t)_{t \geq 0}$ is an increasing family of compact sets in $\overline{\mathbb{H}}$; in addition, $g_t$ is a conformal equivalence of $\mathbb{H} \setminus K_t$ onto $\mathbb{H}$. The families of hulls $(K_t)$ and associated conformal equivalences $(g_t)$ constitute a Loewner chain. If $(W_t/\sqrt{\kappa})$ is a standard Brownian motion (starting from 0), this random Loewner chain defines chordal $\mathrm{SLE}(\kappa)$ in $\mathbb{H}$. It has been proved [20] (see [14] for the case $\kappa = 8$) that there exists a continuous process $(\gamma_t)_{t \geq 0}$ with values in $\overline{\mathbb{H}}$ such that $\mathbb{H} \setminus K_t$ is the unbounded connected component of $\mathbb{H} \setminus \gamma_{[0,t]}$ a.s. This process is the trace of the SLE and it can be recovered from $g_t$ (and therefore from $W_t$) by

$$\gamma_t = \lim_{z \to W_t, z \in \mathbb{H}} g_t^{-1}(z).$$

For any simply connected domain $D$ with two boundary points (or prime ends) $a$ and $b$, chordal SLE$_\kappa$ in $D$ from $a$ to $b$ is defined as $K_t^{(D,a,b)} = h^{-1}(K_t^{(\mathbb{H},0,\infty)})$, where $K_t^{(\mathbb{H},0,\infty)}$ is as above and $h$ is a conformal equivalence of $(D, a, b)$ onto $(\mathbb{H}, 0, \infty)$. This definition is unambiguous up to a linear time change thanks to the scaling property of SLE in the upper half-plane (inherited from the scaling property of the driving process $W_t$).

We now turn to SLE($\kappa, \rho$) processes, defined in [16]. Let $(W_t, O_t)_{t\geq 0}$ be a two-dimensional semimartingale that satisfies the stochastic differential equations (SDEs)

(2.1)
$$dW_t = \sqrt{\kappa}\, dB_t + \frac{\rho}{W_t - O_t}\, dt,$$
$$dO_t = \frac{2}{O_t - W_t}\, dt,$$

where $B$ is a standard Brownian motion and the inequality $W_t \geq O_t$ is valid for all positive times. This process is well defined for $\kappa > 0$, $\rho > -2$. Indeed, we define $Z_t = W_t - O_t$ and note that the process $(Z_t/\sqrt{\kappa})_{t\geq 0}$ must be a Bessel process of dimension $d = 1 + 2(\rho + 2)/\kappa$.

Hence, we can define $Z/\sqrt{\kappa}$ to be such a Bessel process (see, e.g., [19]), then define $O_t = -2\int_0^t du/Z_u$ and finally define $W_t = Z_t + O_t$.

We may therefore define a SLE($\kappa, \rho$) as a stochastic Loewner chain driven by the process $(W_t)$ defined above. The starting point (or rather state) of the process is a couple $(w, o)$ with $w \geq o$, usually set to $(0, 0)$. Then $O_t$ represents the image under the conformal map $g_t$ of the leftmost point of $\partial K_t \cup O_0$. Obviously, for $\rho = 0$, we recover a standard SLE($\kappa$) process.

Later we need left as well as right SLE($\kappa, \rho$) processes. We have just defined left SLE($\kappa, \rho$) processes, which we denote SLE$_l(\kappa, \rho)$ if there is any ambiguity. Right processes are defined in the same fashion except for the condition $W_t \leq O_t$ for all $t \geq 0$; they are denoted SLE$_r(\kappa, \rho)$. Note that left processes starting from $(0, 0)$ are images of the corresponding right processes under the antiholomorphic equivalence $z \mapsto -\bar{z}$.

**3. Hulls and restriction.** In this section we recall some results of [16], which is the basis of this work. Define a + hull as a bounded set $A \subset \mathbb{H}$ such that $A = \overline{A \cap \mathbb{H}}$, $A \cap \mathbb{R} \subset \mathbb{R}_+^*$ and $\mathbb{H} \setminus A$ is (connected and) simply connected (as in [16]). A smooth + hull is a + hull $A$ such that there exists a simple smooth curve $\gamma : [0, 1] \to \mathbb{C}$, $\gamma(0, 1) \subset \mathbb{H}$, $\gamma(0), \gamma(1) \in \mathbb{R}$ and $\mathbb{H} \cap \partial A = \gamma(0, 1)$. If $A$ is a + hull, we denote by $\phi_A$ the conformal equivalence between $\mathbb{H} \setminus A$ and $\mathbb{H}$ that satisfies the hydrodynamic normalization near infinity:

$$\phi_A(z) = z + o(1).$$

Then composition of conformal equivalences gives a semigroup law on hulls:

$$\phi_{A \cdot B} = \phi_B \circ \phi_A.$$

Let $(g_t)$ be a Loewner chain with driving process $(W_t)$ and let $A$ be a hull. If $A \subset g_t^{-1}(\mathbb{H})$ define $h_t = \phi_{g_t(A)}$. Define also $\widetilde{W}_t = h_t(W_t)$. Then $\tilde{g}_t = \phi_{g_t(A)} \circ g_t \circ \phi_A^{-1}$ is itself a time-changed Loewner chain if $t$ is small enough.

Suppose now that the driving process $(W_t)$ of the chain is a semimartingale that satisfies

$$dW_t = \sqrt{\kappa}\, dB_t + b_t\, dt,$$

where $B$ is a standard Brownian motion and $b$ is some bounded progressive process. Obviously, this is applicable to $\mathrm{SLE}(\kappa, \rho)$ processes. Let $z$ be a point in $\mathbb{H} \setminus g_t(A)$ or in a punctured neighborhood of $W_t$ in $\mathbb{R}$. Then the following formulas hold:

$$(3.1) \qquad \partial_t h_t(z) = \frac{2 h'_t(W_t)^2}{h_t(z) - \widetilde{W}_t} - \frac{2 h'_t(z)}{z - W_t},$$

$$(3.2) \qquad \partial_t h'_t(z) = -\frac{2 h'_t(W_t)^2 h'_t(z)}{(h_t(z) - \widetilde{W}_t)^2} + \frac{2 h'_t(z)}{(z - W_t)^2} - \frac{2 h''_t(z)}{z - W_t},$$

$$(3.3) \qquad [\partial_t h_t](W_t) = \lim_{z \to W_t} \left( \frac{2 h'_t(W_t)^2}{h_t(z) - \widetilde{W}_t} - \frac{2 h'_t(z)}{z - W_t} \right) = -3 h''_t(W_t),$$

$$(3.4) \qquad [\partial_t h'_t](W_t) = \lim_{z \to W_t} \partial_t h'_t(z) = \frac{h''_t(W_t)^2}{2 h'_t(W_t)} - \frac{4 h'''_t(W_t)}{3}.$$

Now, using a suitable version of Itô's formula (see [19], Exercise (IV.3.12)), we can derive the SDEs

$$(3.5) \qquad d\widetilde{W}_t = h'_t(W_t)\, dW_t + \left( \frac{\kappa}{2} - 3 \right) h''_t(W_t)\, dt,$$

$$(3.6) \qquad dh'_t(W_t) = h''_t(W_t)\, dW_t + \left( \frac{h''_t(W_t)^2}{2 h'_t(W_t)} + \left( \frac{\kappa}{2} - \frac{4}{3} \right) h'''_t(W_t) \right) dt.$$

Let us recall that the Schwarzian derivative of $h_t$ at $z$ is given by

$$Sh_t(z) = \frac{h'''_t(z)}{h'_t(z)} - \frac{3 h''_t(z)^2}{2 h'_t(z)^2}.$$

Consider now the semimartingale

$$(3.7) \qquad Y_t = h'_t(W_t)^\alpha \exp\left( \lambda \int_0^t \frac{Sh_s(W_s)}{6}\, ds \right).$$

Then Itô's formula yields

$$(3.8) \qquad \frac{dY_t}{Y_t} = \alpha \frac{h''_t(W_t)}{h'_t(W_t)}\, dW_t,$$

where

$$\alpha = \alpha_\kappa = \frac{6-\kappa}{2\kappa},$$

$$\lambda = \lambda_\kappa = \frac{(8-3\kappa)(6-\kappa)}{2\kappa}.$$

We also need results derived in [17] regarding the *Brownian loop soup*. A loop is a continuous map $S^1 \to D$, where $D$ is a simply connected plane domain, and is defined up to reparametrization; the filling $\delta^f$ of a loop $\delta$ is the simply connected compact subset of $D$ that has the same outer boundary as $\delta(S^1)$. The Brownian loop soup is a loop-valued point process parametrized by its intensity $\lambda$. If $A$ is a bounded hull in $\mathbb{H}$ (i.e., $A$ is a compact subset of $\overline{\mathbb{H}}$, $\mathbb{H} \setminus A$ is simply connected and $A = \overline{A \cap \mathbb{H}}$) and $L$ is the random loop soup, we define a random hull $A^L$ as the closure of the complement of the unbounded connected component of $\mathbb{H} \setminus (A \cup \bigcup_{\delta \in L, \delta \cap A \neq \varnothing} \delta)$.

THEOREM 1 ([17]). (i) *Let $(K_t)_{0 \leq t \leq T}$ be a (deterministic) Loewner chain, with driving process $(w_t)$. Let $A$ be a hull in $\mathbb{H}$ not intersecting $K_T$ and let $(h_t)$ be defined as above. If $L$ is a Brownian loop soup in $\mathbb{H}$ with intensity $\lambda$, then*

$$\exp\left(\lambda \int_0^T \frac{Sh_t}{6}(w_t)\, dt\right) = \mathbb{P}(K_T \cap A^L = \varnothing).$$

(ii) *Conformal invariance. Let $L$ be a Brownian loop soup in a domain $D$ with intensity $\lambda$ and let $\phi$ be a conformal equivalence $\phi: D \to D'$. Then $\phi(L)$ has the law of a Brownian loop soup in $D'$ with intensity $\lambda$.*

(iii) *Restriction. Let $L$ be a Brownian soup with intensity $\lambda$ in a domain $D$ and let $A$ be a hull. Then*

$$L' = \{\delta \in L | \delta \cap A = \varnothing\}$$

*has the law of a Brownian soup in $D \setminus A$ with intensity $\lambda$.*

Finally, let us briefly recall from [16] the definition and constructions of one-sided restriction probability measures. For more information, see [16]. For each $\alpha > 0$, there exists exactly one measure on simple curves $\gamma$ from $0$ to $\infty$ in the upper half-plane such that for all + hull $A$,

$$\mathbb{P}(\gamma \cap A = \varnothing) = \phi'_A(0)^\alpha.$$

These are the only measures on curves that satisfy the "one-sided restriction property." The curve $\gamma$ is a sample of the one-sided restriction measure with exponent $\alpha$. For each $\alpha$, various equivalent ways to construct this random curve are described in [16]: First, $\gamma$ is a SLE($8/3, \rho$) process for a well-chosen value of $\rho$. Alternatively, when $\alpha \geq 5/8$, we can add to a SLE($\kappa$) for a well-chosen value of $\kappa$, the set of loops of a Brownian loop soup of intensity $\lambda_\kappa$ that it intersects, and consider the right boundary of the obtained set. We generalize these two constructions herein.

**4. Restriction functionals for SLE($\kappa, \rho$) processes.** The main goal of this section is to derive suitable generalizations of the results in [16] that correspond to the case $\kappa = 8/3$ (in this case $\lambda_\kappa = 0$, the "central charge" is null). Interpretations of this formula in terms of conditioning were discussed in [26].

Throughout this paper, we use the following constants that depend on $\kappa$ and $\rho$:

$$a(\kappa, \rho) = \frac{6 - \kappa}{2\kappa},$$

$$b(\kappa, \rho) = \frac{\rho}{4\kappa}(\rho + 4 - \kappa),$$

$$c(\kappa, \rho) = \frac{\rho}{\kappa},$$

$$\lambda_\kappa = \frac{(8 - 3\kappa)(6 - \kappa)}{2\kappa}.$$

Note that $a$ depends only on $\kappa$.

LEMMA 1. *Suppose that $\kappa > 0$ and $\rho > -2$. Let $(W_t, O_t)$ generate a SLE($\kappa, \rho$) process and let $A$ be a $a$ + hull. Consider the semimartingale*

$$M_t = h'_t(W_t)^a h'_t(O_t)^b \left(\frac{h_t(W_t) - h_t(O_t)}{W_t - O_t}\right)^c \exp\left(\lambda_\kappa \int_0^t \frac{Sh_s}{6}(W_s)\,ds\right).$$

*Then, with the previous choice of constants $a, b, c$ and $\lambda$, the process $(M_t)$, which is well defined up to an a.s. positive stopping time $\tau$, is a local martingale.*

PROOF. This lemma is the natural generalization of [16], Lemma 8.9. The proof is a straightforward application of Itô's formula, which we write down for the sake of completeness. Recall (3.8) and (2.1):

$$\frac{dY_t}{Y_t} = a\frac{h''_t(W_t)}{h'_t(W_t)}\,dW_t$$

$$= a\frac{h''_t(W_t)}{h'_t(W_t)}\left(\sqrt{\kappa}\,dB_t + \frac{\rho}{W_t - O_t}\right)dt.$$

Standard differential calculus yields [see (3.1), (3.2) and (2.1)]

$$dh_t(O_t) = \frac{2h'_t(W_t)^2}{h_t(O_t) - \widetilde{W}_t}\,dt$$

and

$$\frac{dh'_t(O_t)}{h'_t(O_t)} = \left(\frac{2}{(O_t - W_t)^2} - \frac{2h'_t(W_t)^2}{(h_t(O_t) - h_t(W_t))^2}\right) dt.$$

From (3.5) and (2.1), using Itô's formula, we get [we write $U_t = (h_t(W_t) - h_t(O_t))/(W_t - O_t)$]

$$\frac{dU_t}{U_t} = \left(\frac{h'_t(W_t)}{h_t(W_t) - h_t(O_t)} - \frac{1}{W_t - O_t}\right)\sqrt{\kappa}\, dB_t$$
$$+ \left[\frac{(\rho - \kappa)h'_t(W_t)}{(W_t - O_t)(h_t(W_t) - h_t(O_t))} + \frac{\kappa - \rho - 2}{(W_t - O_t)^2}\right.$$
$$\left. + \left(\frac{\kappa}{2} - 3\right)\frac{h''_t(W_t)}{h_t(W_t) - h_t(O_t)} + \frac{2h'_t(W_t)^2}{(h_t(W_t) - h_t(O_t))^2}\right] dt.$$

Then

$$\frac{dM_t}{M_t} = \frac{dY_t}{Y_t} + b\frac{dh'_t(O_t)}{h'_t(O_t)} + c\frac{dU_t}{U_t} + \frac{1}{2}c(c-1)\frac{d\langle U_t\rangle}{U_t^2} + c\frac{d\langle Y_t, U_t\rangle}{Y_t U_t}.$$

By substituting,

$$\frac{dM_t}{M_t} = \left[a\frac{h''_t(W_t)}{h'_t(W_t)} + c\left(\frac{h'_t(W_t)}{h_t(W_t) - h_t(O_t)} - \frac{1}{W_t - O_t}\right)\right]\sqrt{\kappa}\, dB_t$$
$$+ \left[\frac{h''_t(W_t)}{h'_t(W_t)(W_t - O_t)}(a\rho - ac\kappa)\right.$$
$$+ \frac{h''_t(W_t)}{h_t(W_t) - h_t(O_t)}\left(c\left(\frac{\kappa}{2} - 3\right) + ca\kappa\right)$$
$$+ \frac{1}{(W_t - O_t)^2}\left(2b + c(\kappa - \rho - 2) + \frac{1}{2}c(c-1)\kappa\right)$$
$$+ \frac{h'_t(W_t)^2}{(h_t(O_t) - h_t(W_t))^2}\left(-2b + 2c + \frac{1}{2}c(c-1)\kappa\right)$$
$$\left.+ \frac{h'_t(W_t)}{(W_t - O_t)(h_t(W_t) - h_t(O_t))}(c(\rho - \kappa) - c(c-1)\kappa)\right] dt.$$

It is then easy to see that the drift terms vanish for our specific choice of constants $a$, $b$ and $c$. $\square$

To apply the optional stopping theorem, we need two more lemmas. It is convenient to define

$$\alpha(\kappa, \rho) = a(\kappa, \rho) + b(\kappa, \rho) + c(\kappa, \rho) = \frac{(\rho + 2)(\rho + 6 - \kappa)}{4\kappa}.$$

Note that $\alpha(\kappa, 0) = a(\kappa, \rho) = \alpha_\kappa = (6 - \kappa)/2\kappa$. We now assume that $A$ is a smooth + hull and $(M_t)$ is as above.

LEMMA 2. (i) *If $\kappa \leq \frac{8}{3}$, for a left* SLE$(\kappa, \rho)$, *the associated local martingale* $(M_t)$ *is bounded:* $0 \leq M_t \leq 1$.

(ii) *If $\kappa \geq 6$ and $\rho \geq \kappa - 4$, for a right* SLE$(\kappa, \rho)$, *the associated local martingale* $(M_t)$ *is bounded:* $0 \leq M_t \leq 1$.

PROOF. (i) Note that for all $\kappa < 4$ and $\rho > -2$, the exponent $\alpha(\kappa, \rho)$ is positive. Moreover, if $\kappa \leq 8/3$, then $\lambda_\kappa \geq 0$, so that the exponential term is bounded by 1 (since the Schwarzian derivatives are negative in the present case; see, e.g., equation (5.5) in [16]).

Generally speaking, if $B$ is a smooth + hull, and $x$ and $y$ are two real numbers such that $x < y < \inf \mathbb{R} \cap B$, then (see the proof of Lemma 8.10 in [16])

$$1 \geq \phi_B'(x) \geq \phi_B'(y).$$

Since $O_t \leq W_t$ (left SLE), it follows that

$$h_t'(O_t) \geq \frac{h_t(W_t) - h_t(O_t)}{W_t - O_t} \geq h_t'(W_t).$$

We split the proof into different cases, according to the signs of the constants $b$ and $c$:

- Suppose first that $b(\kappa, \rho) \leq 0$. Then, also, $c(\kappa, \rho) \leq 0$. Recall that $a + b + c > 0$. Then,

$$M_t \leq h_t'(W_t)^a h_t'(W_t)^b h_t'(W_t)^c \leq h_t'(W_t)^{\alpha(\kappa, \rho)} \leq 1.$$

- Suppose now that $b(\kappa, \rho) > 0$ and $c(\kappa, \rho) \geq 0$. Then, trivially,

$$M_t \leq h_t'(W_t)^a \leq 1.$$

- Suppose finally that $b(\kappa, \rho) > 0$ and $c(\kappa, \rho) < 0$. The hull $g_t(A)$ is a smooth hull, so it has a Loewner parametrization, that is, there is a continuous real-valued function $(x_s)_{0 \leq s \leq S}$ such that $g_t(A) = \widetilde{K}_S$, the Loewner hull associated with $x$ at time $S$. Let $(\tilde{g}_s)_{0 \leq s \leq S}$ be the corresponding conformal equivalences. Then, if $o_s = \tilde{g}_s(O_t)$ and $w_s = \tilde{g}_s(W_t)$, we get (see [16], proof of Theorem 8.4)

$$M_t = \exp\left(-2 \int_0^S \left(\frac{a(\kappa, \rho)}{(x_s - w_s)^2} + \frac{c(\kappa, \rho)}{(x_s - w_s)(x_s - o_s)} + \frac{b(\kappa, \rho)}{(x_s - o_s)^2}\right) ds\right)$$
$$\times \exp\left(\lambda_\kappa \int_0^S \frac{Sh_u(W_u)}{6} du\right).$$

Let $y = (x_s - w_s)/(x_s - o_s)$. Then it is sufficient to prove that

$$b(\kappa, \rho) y^2 + c(\kappa, \rho) y + a(\kappa, \rho) \geq 0 \qquad \forall y \in [0, 1],$$

but in this case, $-c/(2b) \geq 1$ and, hence,

$$\min_{y \in [0,1]} (a + cy + by^2) = a + b + c = \alpha(\kappa, \rho) > 0.$$

We can conclude that $0 \leq M_t \leq 1$. A slight modification of the argument gives, in fact, that for some positive $\varepsilon$,

$$0 \leq M_t \leq h'(W_t)^\varepsilon.$$

(ii) In the case where $\kappa \geq 6$ and $\rho \geq \kappa - 4$, $\alpha(\kappa, \rho)$ and $c(\kappa, \rho)$ are positive, $a(\kappa, \rho)$ is nonpositive and $b(\kappa, \rho)$ is nonnegative. We are now dealing with a right SLE process, $W_t \leq O_t$; hence,

$$h'_t(O_t) \leq \frac{h_t(W_t) - h_t(O_t)}{W_t - O_t} \leq h'_t(W_t).$$

Again, $\lambda_\kappa$ is positive, so that the exponential factor is bounded by 1. It readily follows that

$$M_t \leq \left( \frac{h_t(W_t) - h_t(O_t)}{W_t - O_t} \right)^{a+b+c} \leq h'_t(W_t)^{\alpha(\kappa, \rho)}.$$

This concludes the proof. Note that we have, in fact, proved in all these cases the existence of a positive $\varepsilon$ such that $M_t \leq h'_t(W_t)^\varepsilon$. $\square$

Recall that $\gamma$ denotes the trace of a SLE. Somewhat loosely, we also use $\gamma$ to designate the closed set $\gamma_{[0,\infty[}$. For a smooth + hull $A$, we define a bounded martingale $(M_t)$ as in Lemma 1.

LEMMA 3. (i) *If $\kappa \leq \frac{8}{3}$, for a left* SLE($\kappa, \rho$), *the martingale $(M_t)$ converges a.s.:*

$$M_t \to \mathbf{1}_{\gamma \cap A = \varnothing} \exp\left( \lambda_\kappa \int_0^\infty \frac{Sh_s(W_s)}{6} ds \right).$$

(ii) *If $\kappa \geq 6$, $\rho = \kappa - 4$, for a right* SLE($\kappa, \rho$), *the martingale $(M_t)$ converges a.s.:*

$$M_t \to \mathbf{1}_{\gamma \cap A = \varnothing} \exp\left( \lambda_\kappa \int_0^\infty \frac{Sh_s(W_s)}{6} ds \right).$$

PROOF. (i) In this case, we have proved that $M_t \leq h'_t(W_t)^\varepsilon$ for some $\varepsilon > 0$. From [16], Lemma 8.3, the trace $\gamma$ a.s. does not intersect $(0, \infty)$. Then we can apply Lemmas 6.2 and 6.3 of [16] to get the result.

(ii) If $\kappa \geq 6$ and $\rho = \kappa - 4$, we have seen that $0 \leq M_t \leq h'_t(W_t)^{\alpha(\kappa, \rho)}$. Moreover, $(O_t - W_t)$ is a transient Bessel process of dimension $(3 - 4/\kappa)$, so the trace $\gamma$ a.s. does not intersect $(0, \infty)$. Then Lemma 6.3 of [16] tells us that on the event $\{\gamma \cap A \neq \varnothing\}$, $M_t \to 0$ as $t \nearrow \tau_A$, the first time for which the trace encounters the hull. On the event $\{\gamma \cap A = \varnothing\}$ or $\{\tau_A = \infty\}$, as before, $h'_t(W_t) \to 1$ as $t$ goes to infinity. Note that $b(\kappa, \kappa - 4) = 0$, so we have to prove that

$$\frac{h_t(W_t) - h_t(O_t)}{W_t - O_t} \xrightarrow[t \to \infty]{} 1. \tag{4.1}$$

Here we need to adapt the proof of Lemma 6.2 in [16]. Heuristically, seen from $A$, $W_t$ goes to infinity while $O_t$ stays bounded [in particular, $h'_t(O_t)$ tends to a nondegenerate limit, which is why we do not consider the case $\rho > \kappa - 4$]. Let $A_t = g_t(A)$, $Z_t$ be the leftmost point of $A_t$ and let $d_t = \inf(r, A_t \subset D(Z_t, r))$. We also define $D_t = \{z \in \overline{\mathbb{H}}, |z - Z_t| \leq d_t\}$ and $O'_t = \min(O_t, Z_t - d_t)$. The extremal distance between $g_t^{-1}((-\infty, W_t))$ and $\partial A$ in $\mathbb{H} \setminus (K_t \cup A)$ goes to infinity, while the extremal distance between $g_t^{-1}((-\infty, O_t))$ and $\partial A$ stays bounded. Indeed, since we are dealing with a right SLE, $O_t$ is the right image of 0 under $g_t$. This can be translated into

$$\frac{d_t}{Z_t - W_t} \to 0, \qquad \liminf \frac{d_t}{Z_t - O_t} > 0,$$

which implies that $d_t/(O_t - W_t) \to 0$. We have already seen that

$$\frac{h_t(O_t) - h_t(W_t)}{O_t - W_t} \leq 1.$$

However, $\phi_{D_t} = \phi_{h_t(D_t)} \circ h_t$, so

$$\frac{\phi_{D_t}(O'_t) - \phi_{D_t}(W_t)}{h_t(O'_t) - h_t(W_t)} \leq 1.$$

Now, by scaling, $\phi_{D_t}(z) = d_t \phi_D((z - Z_t)/d_t) + Z_t$ for a fixed function $\phi_D$ which satisfies the hydrodynamic normalization at infinity: $\phi_D(z) = z + o(1)$. Hence,

$$O_t - W_t \geq h_t(O_t) - h_t(W_t)$$
$$\geq \phi_{D_t}(O'_t) - \phi_{D_t}(W_t) = d_t\left(\frac{O_t - W_t}{d_t} + O(1)\right),$$

so we can conclude that (4.1) holds. □

From the previous lemmas, we get immediately:

PROPOSITION 1. (i) *Consider a* $\mathrm{SLE}_l(\kappa, \rho)$, $\kappa \leq \frac{8}{3}$, *and a smooth + hull* $A$. *Then*

$$\phi'_A(0)^{\alpha(\kappa,\rho)} = \mathbb{E}\bigg(\mathbf{1}_{\gamma \cap A = \varnothing} \exp\bigg(\lambda_\kappa \int_0^\infty \frac{Sh_s(W_s)}{6} ds\bigg)\bigg).$$

(ii) *Consider a* $\mathrm{SLE}_r(\kappa, \kappa - 4)$, $\kappa \geq 6$, *and a smooth + hull* $A$. *Then*

$$\phi'_A(0)^{1/2 - 1/\kappa} = \mathbb{E}\bigg(\mathbf{1}_{\gamma \cap A = \varnothing} \exp\bigg(\lambda_\kappa \int_0^\infty \frac{Sh_s(W_s)}{6} ds\bigg)\bigg).$$

In the case $\kappa \geq 6$ and $\rho > \kappa - 4$, we can derive a formula that involves the nondegenerate value of $h'_t(O_t)$; for more on this topic, see [26].

The first statement shows that for all $\alpha > 0$, we can construct a sample of the one-sided restriction measure by adding to a SLE($\kappa, \rho$) a Poisson cloud of Brownian bubbles of intensity $\lambda_\kappa$, when $\kappa \leq 8/3$ and $\alpha = \alpha(\kappa, \rho)$.

We now enunciate a corollary that provides important support for a duality conjecture. Consider a SLE$_r(\kappa, \rho)$ with $\kappa \geq 6$ and $\rho \geq \kappa - 4$. Then it is easy to see from the previous results that for any + hull $A$, $K_\infty \cap A = \varnothing$ with positive probability. Thus it makes sense to define the right boundary of $K_\infty$. The purpose of duality is to identify this boundary as a process.

COROLLARY 1. *Let $\kappa \geq 6$ and $\kappa' = 16/\kappa$. Let $\delta$ be the right boundary of a* SLE$_r(\kappa, \kappa - 4)$ *and let $\gamma'$ be the trace of a* SLE$_l(\kappa', \frac{\kappa'-4}{2})$. *Let $L$ be an independent Brownian loop soup in $\mathbb{H}$ with intensity $\lambda_\kappa$. Then for any + hull $A$,*

$$\mathbb{P}(\delta \cap A^L = \varnothing) = \mathbb{P}(\gamma' \cap A^L = \varnothing) = \phi'_A(0)^{1/2 - 1/\kappa}.$$

PROOF. The result is immediate using the properties of the loop soup. Indeed, if $l$ is any loop in $\mathbb{H}$ that intersects $A$, for obvious topological reasons,

$$\{l \cap \delta \neq \varnothing\} = \{l \cap K_\infty \neq \varnothing\}.$$

Moreover, $\delta$ (resp. $K_\infty$) intersects $A^L$ if and only if it intersects a loop $l \in L$. Then for a smooth + hull $A$, almost surely

$$\exp\left(\lambda_\kappa \int_0^\infty \frac{Sh_s(w_s)}{6} ds\right) = \exp\left(\lambda_\kappa \int_0^\infty \frac{S\tilde{h}_s(\tilde{w}_s)}{6} ds\right),$$

where the left-hand side corresponds to the SLE$_r(\kappa, \kappa - 4)$ process itself, while the right-hand side corresponds to its right boundary process with Loewner parametrization. This implies the result for general + hulls (see Lemma 2.1 of [16]). □

This suggests the following conjecture (with the same notation):

CONJECTURE 2. *The simple curves $\delta$ and $\gamma'$ have the same law.*

This conjecture should also hold for $\kappa \in (4, 6)$. Note that when $\kappa = 4$, it trivially holds.

**5. Some properties of SLE($\kappa, \kappa - 4$) processes.** We have just seen a precise duality conjecture that involves SLE($\kappa, \kappa - 4$) processes. We now study these processes, which satisfy particular properties. In some sense, we may see what follows as a rephrasing of the following well-known properties of Bessel processes: When $d < 2$, a Bessel process conditioned never to hit the origin has the law of a Bessel process of dimension $4 - d$. Here, this is translated into the fact that $\text{SLE}_r(\kappa, \kappa - 4)$ is a $\text{SLE}(\kappa)$ conditioned not to hit the positive half-line.

PROPOSITION 2. *Let $\kappa > 4$. A $\text{SLE}(\kappa)$ conditioned not to absorb $x > 0$ has the law of a $\text{SLE}_r(\kappa, \kappa - 4)$ starting from $(0, x)$.*

This conditioning with respect to a zero probability event holds only under an appropriate limiting procedure (see the proof).

PROOF OF PROPOSITION 2. Let $(W_t)$ be the driving process of the $\text{SLE}(\kappa)$, let $dW_t = \sqrt{\kappa}\, dB_t$, and let $L$ be a (large) negative number. Let $L_t = g_t(L)$ and $x_t = g_t(x)$; these processes are defined up to time $\tau$, at which either $L$ or $x$ is swallowed. Let $h$ be a function on $(0, 1)$ that satisfies the ODE

$$\frac{\kappa}{4} h''(z) + \left(\frac{1}{z} - \frac{1}{1-z}\right) h'(z) = 0.$$

Then, according to [25], Proposition 3.3, $h((W_t - L_t)/(x_t - L_t))$ is a (local) martingale. Now let $h(1) = 0$, $h(0) = 1$ and $Z_t = (W_t - L_t)/(x_t - L_t)$. The following SDE holds:

$$dZ_t = \frac{\sqrt{\kappa}\, dB_t}{x_t - L_t} + \frac{2\,dt}{(x_t - L_t)^2}\left(\frac{1}{Z_t} - \frac{1}{1 - Z_t}\right).$$

Suppose that the processes $B., W., L.$ and $x.$ are defined on a filtered probability space $(\Omega, \mathcal{F}_t, \mathbf{P})$, so $B$ is a standard Brownian motion under $\mathbf{P}$, and let $\mathbf{Q}$ be the conditional measure $\mathbf{Q} = \mathbf{P}(\cdot | L \text{ is swallowed before } x)$. Then

$$\frac{d\mathbf{Q}}{d\mathbf{P}}\bigg|_{\mathcal{F}_t} = \frac{h(Z_t)\mathbf{1}_{\tau \geq t} + \mathbf{1}_{\tau \leq t, \gamma_\tau = L}}{h(-L/(x-L))} = D_t.$$

According to Girsanov's theorem, if $\widetilde{B}$ satisfies $d\widetilde{B}_t = dB_t - D_t^{-1} d\langle B, D\rangle_t$, then under $\mathbf{Q}$, $\widetilde{B}$ is a continuous local martingale with quadratic variation, $\langle \widetilde{B}\rangle_t = \langle B\rangle_t = t$, so $\widetilde{B}$ is a standard Brownian motion under $\mathbf{Q}$. So we have

$$dB_t = d\widetilde{B}_t + \frac{\sqrt{\kappa}}{x_t - L_t} \frac{h'}{h}(Z_t)\, dt.$$

As $L$ goes to $-\infty$, $Z_t$ converges to 1. It is easily seen that as $z \to 1$,

$$\frac{h'(z)}{h(z)} \sim -\frac{1 - 4/\kappa}{1 - z}.$$

Indeed $h'(z) = c(z(1-z))^{-4/\kappa}$, so that $h'(z) \sim c(1-z)^{-4/\kappa}$ and $h(z) \sim -(c/(1-4/\kappa))(1-z)^{1-4/\kappa}$ as $z \to 1$. So when $L \to -\infty$, which corresponds to conditioning by the event of zero probability $\{K_\infty \cap (x, \infty) = \varnothing\}$, we get the SDE

$$dW_t = \sqrt{\kappa}\, d\widetilde{B}_t + \frac{\kappa - 4}{W_t - x_t}\, dt$$

and $x_t$ satisfies the Loewner equation, which is the definition of a SLE($\kappa, \kappa - 4$) process under $\mathbf{Q}$.  □

PROPOSITION 3. *Let $\kappa > 4$. A SLE($\kappa, \kappa - 4$) starting from $(0, x)$ and conditioned not to absorb $y$, $0 < y < x$, has the law of a SLE$_r$($\kappa, \kappa - 4$) starting from $(0, y)$.*

PROOF. Making use of the previous result, we can formally interpret this result as

$$\begin{aligned}
\text{SLE}_{(0,y)}(\kappa, \kappa - 4) &= (\text{SLE}(\kappa)|y \text{ is not absorbed}) \\
&= (\text{SLE}(\kappa)|x \text{ and } y \text{ are not absorbed}) \\
&= ((\text{SLE}(\kappa)|x \text{ is not absorbed})|y \text{ is not absorbed}) \\
&= (\text{SLE}_{(0,x)}(\kappa, \kappa - 4)|y \text{ is not absorbed}).
\end{aligned}$$

We can derive a proof for this fact along the lines of the previous proposition, that is, using Girsanov's theorem. It is easy to check that if $(W_t, O_t)$ is the driving process of a SLE($\kappa, \kappa - 4$) starting from $(0, x)$, and $0 < y < x$, then

$$\left( \frac{g_t(y) - W_t}{g_t(x) - W_t} \right)^{1 - 4/\kappa}$$

is a bounded martingale. Consequently, if $dW_t = \sqrt{\kappa}\, dB_t + \frac{\kappa - 4}{W_t - g_t(x)}\, dt$, where $B$ is a standard Brownian motion defined on the filtered probability space $(\Omega, \mathcal{F}_t, \mathbf{P})$, and if $\mathbf{Q}$ denotes the conditional measure: $\mathbf{Q} = \mathbf{P}(\cdot | y \text{ is not swallowed})$, then

$$\frac{d\mathbf{Q}}{d\mathbf{P}}_{|\mathcal{F}_t} = c \left( \frac{g_t(y) - W_t}{g_t(x) - W_t} \right)^{1 - 4/\kappa} = D_t,$$

where $c = (x/y)^{1 - 4/\kappa}$. If $\widetilde{B}$ satisfies $d\widetilde{B}_t = dB_t - D_t^{-1} d\langle B, D \rangle_t$, then $\widetilde{B}$ is a standard Brownian motion under $\mathbf{Q}$. We can compute

$$dD_t = c \left( 1 - \frac{4}{\kappa} \right) \left( \frac{g_t(y) - W_t}{g_t(x) - W_t} \right)^{-4/\kappa} \frac{g_t(y) - g_t(x)}{(g_t(x) - W_t)^2} \sqrt{\kappa}\, dB_t$$

so that

$$dW_t = \sqrt{\kappa}\, dB_t + \frac{\kappa - 4}{W_t - g_t(x)}\, dt$$

$$= \sqrt{\kappa}\, d\widetilde{B}_t + \sqrt{\kappa}\frac{d\langle B, D\rangle_t}{D_t} + \frac{\kappa - 4}{W_t - g_t(x)}\, dt$$

$$= \sqrt{\kappa}\, d\widetilde{B}_t + (\kappa - 4)\left(\frac{1}{W_t - g_t(x)} + \frac{g_t(y) - g_t(x)}{(g_t(x) - W_t)(g_t(y) - W_t)}\right) dt$$

$$= \sqrt{\kappa}\, d\widetilde{B}_t + \frac{\kappa - 4}{W_t - g_t(y)}\, dt.$$

Under **Q**, this defines a SLE$(\kappa, \kappa - 4)$ starting from $(0, y)$, which concludes the proof. $\square$

Consider now a SLE$_r(\kappa, \kappa - 4)$ starting from $(0, 0^+)$ with trace $\gamma$, right boundary $\delta$ and driving process $(W_t, O_t)$. Let $(\mathcal{F}_t)$ be the natural filtration of the Brownian motion $(B_t)$ that drives the SDE of $(W_t)$. Now $\delta$ is a simple curve that can be parametrized so that $\operatorname{cap}(\delta_{[0,u]}) = 2u$, where cap is the half-space capacity seen from infinity (in the terminology of [11]). Let $\tau_u$ be the first time at which the portion $\delta_{[0,u]}$ of the boundary is completed; obviously this is not a stopping time. Formally, we can define $(\mathcal{D}_u)_{u \geq 0}$ as the filtration generated by $(\delta_u)_{u \geq 0}$, as well as a finer filtration $(\widetilde{\mathcal{D}}_u)_{u \geq 0} = (\mathcal{F}_{\tau_u})_{u \geq 0}$.

Let $u > 0$ be fixed. Consider a time $t > 0$. A SLE$_r(\kappa, \kappa - 4)$ starting from $(0, 0^+)$ is the concatenation of the hull $K_t$ with the hull produced by an independent SLE$_r(\kappa, \kappa - 4)$ starting from $(W_t, O_t)$. Then $\tau_u \leq t$ if and only if the right boundary of $K_t$ [i.e., $g_t^{-1}([W_t, O_t])$] has capacity larger than $2u$ and if the future hull does not swallow $g_t(\delta_u)$. So conditionally on $(W_t, O_t, \tau_u \leq t)$, the future is a SLE$_r(\kappa, \kappa - 4)$ starting from $(W_t, g_t(\delta_u))$, independent from the past. Then $(g_{\tau_u}(K_{\tau_u + t}))_{t \geq 0}$ has the law of a SLE$_r(\kappa, \kappa - 4)$ starting from $(W_{\tau_u}, W_{\tau_u}^+)$ and is independent from $\widetilde{\mathcal{D}}_u$ conditionally on $W_{\tau_u}$. Denote by $K'$ the closed subset of $\mathbb{H}$ swallowed by this process, $K' = \overline{g_{\tau_u}(K_\infty \setminus K_{\tau_u})}$, so that $K_\infty = \bigcup_{t \geq 0} K_t$ is the concatenation of $K_{\tau_u}$ and $K'$, and $(K' - W_{\tau_u})$ is a copy of $K_\infty$ independent from $\mathcal{F}_{\tau_u}$. Thus $\tau_u$ can be called a regeneration time (it is a splitting time in Williams' terminology). To get a full path decomposition, we need to describe $K_{\tau_u}$.

We have conjectured that $(\delta_v)_{0 \leq v \leq u}$ is a SLE$_l(\kappa', \frac{\kappa'-4}{2})$ starting from $(0, 0^-)$ stopped at a fixed time $u$. Let $(\widetilde{W}_v, \widetilde{O}_v)_v$ be its driving process and let $(\tilde{g}_v)_{0 \leq v \leq u}$ be the associated conformal equivalences. Now we can write $K_{\tau_u} = \delta_{[0,u]} \cdot H_u$, where $H_u = \tilde{g}_u^{-1}(K_{\tau_u} \setminus \delta_{[0,u]})$. We use the (loose) notation $\delta_u = \delta_{[0,u]}$. Invoking "conformal invariance," we can conjecture that $H_u$ is independent from $\mathcal{D}_u$ conditionally on $(\widetilde{W}_u, \widetilde{O}_u)$. Then, if $A$ is a smooth + hull, $L$ and $L'$ are independent loop soups in $\mathbb{H}$ with intensity $\lambda_\kappa$, and as usual $h_t = \phi_{g_t(A)}$ and $\tilde{h}_u = \phi_{\tilde{g}_u(A)}$, using the restriction property of the loop soup and restriction formulas for SLE$_r(\kappa, \kappa - 4)$ and SLE$_l(\kappa', \kappa'/2 - 2)$, we

get

$$\phi'_A(0)^{1/2-1/\kappa} = \mathbb{P}((\delta_u \cdot H_u \cdot K') \cap A^L = \varnothing) = \mathbb{E}(\mathbf{1}_{\delta_u \cdot H_u \cap A^L = \varnothing} \mathbf{1}_{K' \cap g_{\tau_u}(A)^{L'}})$$

$$= \mathbb{E}(\mathbf{1}_{(\delta_u \cdot H_u \cap A^L) = \varnothing} h'_{\tau_u}(W_{\tau_u})^{1/2-1/\kappa})$$

$$= \mathbb{E}(\mathbf{1}_{\delta_u \cap A^L = \varnothing} \mathbf{1}_{H_u \cap \tilde{g}_u(A)^{L'} = \varnothing} h'_{\tau_u}(W_{\tau_u})^{1/2-1/\kappa})$$

and

$$\phi'_A(0)^{1/2-1/\kappa} = \mathbb{P}(\delta_u \cdot \tilde{g}_u(\delta_{[u,\infty]}) \cap A^L = \varnothing)$$

$$= \mathbb{E}(\mathbf{1}_{\delta_u \cap A^L = \varnothing} \mathbf{1}_{\tilde{g}_u(\delta_{[u,\infty]}) \cap \tilde{g}_u(A)^{L'} = \varnothing})$$

$$= \mathbb{E}\left(\mathbf{1}_{\delta_u \cap A^L = \varnothing} \tilde{h}'_u(\widetilde{W}_u)^{a_{\kappa'}} \tilde{h}'_u(\widetilde{O}_u)^{-(\kappa'-4)^2/(16\kappa')} \right.$$

$$\left. \times \left(\frac{\tilde{h}_u(\widetilde{W}_u) - \tilde{h}_u(\widetilde{O}_u)}{\widetilde{W}_u - \widetilde{O}_u}\right)^{(\kappa'-4)/(2\kappa')}\right).$$

This computation leads us to conjecture that if $B$ is any $+$ hull [in particular, $B = \tilde{g}_u(A)$], then, conditionally on $(\widetilde{W}_u, \widetilde{O}_u) = (w, o)$, $H_u$ satisfies

$$\phi'_B(w)^{a_{\kappa'}} \phi'_B(o)^{-(\kappa'-4)^2/(16\kappa')} \left(\frac{\phi_B(w) - \phi_B(o)}{w - o}\right)^{(\kappa'-4)/(2\kappa')}$$

$$= \mathbb{E}(\mathbf{1}_{H_u \cap B^L = \varnothing} h'(\phi_{H_u}(w^+))^{1/2-1/\kappa}),$$

where $h = \phi_{\phi_{H_u}(B)}$.

In the next section we define random hulls that satisfy this particular restriction formula.

**6. Generalized SLE($\kappa, \underline{\rho}$) processes.** Let $\kappa > 0$ and let $\underline{\rho}$ be a multi-index, that is,

$$\underline{\rho} \in \bigcup_{i \geq 0} \mathbb{R}^i.$$

Let $k$ be the length of $\underline{\rho}$; if $k = 0$, we simply define SLE($\kappa, \varnothing$) as a standard SLE($\kappa$). If $k > 0$, suppose the existence of processes $(W_t)_{t \geq 0}$ and $(Z_t^{(i)})_{t \geq 0}$, $i \in \{1, \ldots, k\}$, that satisfy the SDEs

(6.1)
$$dW_t = \sqrt{\kappa} \, dB_t + \sum_{i=1}^{k} \frac{\rho_i}{W_t - Z_t^{(i)}} \, dt,$$

$$dZ_t^{(i)} = \frac{2}{Z_t^{(i)} - W_t} \, dt$$

and such that the processes $(W_t - Z_t^{(i)})$ do not change sign. Then we define the $\mathrm{SLE}(\kappa, \underline{\rho})$ process that starts from $(w, z_1, \ldots, z_k)$ as a Schramm–Loewner evolution, the driving process of which has the same law as $(W_t)$ defined above, with $W_0 = 0$ and $Z_0^{(i)} = z_i$. Obviously, for $k = 1$, $\underline{\rho} = (\rho)$, we recover the definition of a $\mathrm{SLE}(\kappa, \rho)$ process.

LEMMA 4. *Let $\kappa > 0$ and $\underline{\rho} = (\rho_1, \ldots, \rho_k)$. Suppose that the $\mathrm{SLE}(\kappa, \underline{\rho})$ process exists up to time $\tau$ and let*

$$(W_t, Z_t^{(1)}, \ldots, Z_t^{(k)})$$

*be its driving mechanism. Let $A$ be a $+$ hull and let $(h_t)$ be the associated family of conformal equivalences. For $i \in \{1, \ldots, k\}$, $t < \tau$, define*

$$M_t^{(i)} = h_t'(Z_t^{(i)})^{\rho_i(\rho_i + 4 - \kappa)/(4\kappa)} \left( \frac{h_t(W_t) - h_t(Z_t^{(i)})}{W_t - Z_t^{(i)}} \right)^{\rho_i/\kappa}.$$

*For $i, j \in \{1, \ldots, k\}$, $i < j$, define also*

$$M_t^{(i,j)} = \left( \frac{h_t(Z_t^{(i)}) - h_t(Z_t^{(j)})}{Z_t^{(i)} - Z_t^{(j)}} \right)^{\rho_i \rho_j / (2\kappa)}.$$

*Finally, let*

$$M_t^{\varnothing} = h_t'(W_t)^{\alpha_\kappa} \exp\left( \lambda_\kappa \int_0^t \frac{Sh_s(W_s)}{6} \, ds \right).$$

*Then the semimartingale*

$$M_t = M_t^{\varnothing} \prod_{1 \leq i \leq k} M_t^{(i)} \prod_{1 \leq i < j \leq k} M_t^{(i,j)}$$

*is a local martingale. Moreover, the sum of all exponents in this local martingale equals $\alpha(\kappa, \rho_1 + \cdots + \rho_k)$.*

PROOF. This generalization of Lemma 1 relies on the results recalled in Section 3. First, $M^{\varnothing}$ is the semimartingale formerly denoted by $Y$, so

$$\frac{dM_t^{\varnothing}}{M_t^{\varnothing}} = \alpha_\kappa \frac{h_t''(W_t)}{h_t'(W_t)} \, dW_t.$$

Standard differential calculus yields

$$\frac{dM_t^{(i,j)}}{M_t^{(i,j)}} = \frac{\rho_i \rho_j}{\kappa} \left( \frac{1}{(Z_t^{(i)} - W_t)(Z_t^{(j)} - W_t)} - \frac{h_t'(W_t)^2}{(h_t(Z_t^{(i)}) - \widetilde{W}_t)(h_t(Z_t^{(j)}) - \widetilde{W}_t)} \right) dt.$$

Applying Itô's formula, we get

$$\frac{dM_t^{(i)}}{M_t^{(i)}} = \frac{\rho_i}{\kappa} \left( \frac{h_t'(W_t)}{h_t(W_t) - h_t(Z_t^{(i)})} - \frac{1}{W_t - Z_t^{(i)}} \right) \sqrt{\kappa}\, dB_t$$
$$+ \frac{\rho_i}{\kappa} \bigg[ \left( \frac{\kappa}{2} - 3 \right) \frac{h_t''(W_t)}{h_t(W_t) - h_t(Z_t^{(i)})}$$
$$+ \sum_{j \neq i} \rho_j \bigg( \frac{h_t'(W_t)}{(h_t(W_t) - h_t(Z_t^{(i)}))(W_t - Z_t^{(j)})}$$
$$- \frac{1}{(W_t - Z_t^{(i)})(W_t - Z_t^{(j)})} \bigg) \bigg] dt.$$

Since the "rectangular" semimartingales $M^{(i,j)}$ have no quadratic variation, we get

$$\frac{dM_t}{M_t} = \frac{dM_t^\varnothing}{M_t^\varnothing} + \sum_i \frac{dM_t^{(i)}}{M_t^{(i)}} + \sum_{i<j} \frac{dM_t^{(i,j)}}{M_t^{(i,j)}}$$
$$+ \sum_i \frac{d\langle M_t^\varnothing, M_t^{(i)} \rangle}{M_t^\varnothing M_t^{(i)}} + \sum_{i<j} \frac{d\langle M_t^{(i)}, M_t^{(j)} \rangle}{M_t^{(i)} M_t^{(j)}}.$$

There remains only to check that all the drift terms cancel out. □

The reader with a liking for generality will see that the lemma is a statement on some cancellations for certain quadratic forms; hence, it formally holds for more general SLE processes parametrized by a signed measure $\mu$ (the case we consider corresponds to $\mu = \sum \rho_i \delta_{z_i}$). The main difficulty is establishing that the computations actually make sense.

Let us mention a situation where such processes arise naturally. Consider a chordal SLE($\kappa$) in the half-plane $\mathbb{H}$; let $x < 0 < y$. Then, conditioning the SLE never to swallow either $x$ or $y$, we get a SLE($\kappa, \kappa - 4, \kappa - 4$) starting from $(0, x, y)$. This is a consequence of [13], Theorem 3.1. Recall the discussion of bilateral restriction measures in [16]. It is then straightforward to extend proofs of Lemmas 2 and 3 to get a bilateral analogue to Proposition 1:

COROLLARY 2. *Let $K$ be the hull generated by a SLE($\kappa, \kappa - 4, \kappa - 4$) starting from $(0, 0^-, 0^+)$ with $\kappa \geq 6$. Let $L$ be an independent loop soup with intensity $\lambda_\kappa$. Then the law of $K^L$ is the bilateral restriction measure with exponent $\alpha(\kappa, 2\kappa - 8) = (\kappa - 3)(\kappa - 2)/(2\kappa)$. In particular, for $\kappa = 6$, $K$ has the law of the filling of a Brownian excursion.*

PROOF. By construction, the law of $K^L$ is invariant under $z \mapsto -\bar{z}$. So, by Proposition 3.3 of [16], we have only to check that

$$\mathbb{P}(K^L \cap A = \varnothing) = \phi'_A(0)^{\alpha(\kappa, 2\kappa - 8)}$$

for all smooth + hulls $A$. From the previous proposition, $(M_t)$ is a local martingale. As $t \searrow 0$, $M_t \to M_0 = \phi'_A(0)^{\alpha(\kappa, 2\kappa - 8)}$. Adapting the proofs of Lemmas 2 and 3, we get that $(M_t)$ is bounded and has a.s. limit

$$\mathbf{1}_{K \cap A = \varnothing} \exp\left(\lambda_\kappa \int_0^\infty \frac{Sh_s(W_s)}{6} ds\right).$$

Given the interpretation of the Schwarzian integral in terms of loop soup, this entails the restriction formula. For $\kappa = 6$, $\alpha(\kappa, 2\kappa - 8) = 1$ and it is known that the filling of the path of the Brownian excursion in $\mathbb{H}$ is the bilateral restriction measure with exponent 1 ([16], Proposition 4.1). □

See also [24] for a geometric excursion theory for the Brownian excursion. Recall the discussion at the end of the previous section. We can now make the following statement:

PROPOSITION 4. *Let $(K_t)_{t \geq 0}$ be a $\mathrm{SLE}_r(\kappa, \kappa - 4)$ starting from $(0, 0^+)$, with driving process $(W_u, O_u)$, and let $\tau_u$ be a corresponding regeneration time. Let $(\delta_v)_{v \geq 0}$ be a $\mathrm{SLE}_l(\kappa', \frac{\kappa' - 4}{2})$ starting from $(0, 0^-)$, with driving process $(W'_v, O'_v)_{v \geq 0}$. Let $(H_s)$ be a $\mathrm{SLE}(\kappa, \frac{\kappa}{2} - 4, -\frac{\kappa}{2})$, independent of the former conditionally on $(W'_u, O'_u)$, starting from $(O'_u, (O'_u)^+, W'_u)$. Denote by $\sigma$ the first time at which $W'_u$ is swallowed [the chain $(H_s)$ is only defined up to time $\sigma$]. Finally, let $K_1 = K_{\tau_u}, w_1 = W_{\tau_u}$ and $K_2 = \delta_u \cdot H_\sigma, w_2 = \phi_{K_2}(W'_u)$. Then, for $j = 1, 2$:*

(i) *The capacity of the right boundary of $K_j$ is $2u$ a.s.*

(ii) *For any smooth + hull $A$, if $\alpha = \frac{1}{2} - \frac{1}{\kappa}$ and $L$ is an independent loop soup with intensity $\lambda_\kappa$,*

$$\phi'_A(0)^\alpha = \mathbb{E}(\phi'_{\phi_{K_j}(A)}(w_j)^\alpha \mathbf{1}_{K_j \cap A^L = \varnothing}).$$

This supports the following speculation:

CONJECTURE 3. *The random hulls $K_1$ and $K_2$ are identical in law (at least in the Hausdorff sense).*

This last conjecture is more precise than the general duality conjecture for SLE.

Finally, let us discuss some ideas that will be developed in [7]. We have seen that $\mathrm{SLE}(\kappa, \kappa - 4)$ processes appear quite naturally in the study of restriction formulas. Although their first properties, as presented in Section 5, are independent of the restriction framework (and are valid in the whole range $\kappa > 4$), the study of $\text{SLE}_r(\kappa, \kappa - 4)$ can be carried a bit further. In particular, if $\sigma_t$ denotes the first time after $t$ spent by the trace on the final right boundary ($\delta$ in the notation of Section 5; $\sigma_t$ is obviously not a stopping time), then we can prove that $X_t = g_t(\gamma_{\sigma_t}) - W_t$ defines a Bessel process of dimension $(1 - 4/\kappa)$. This process vanishes exactly when the trace $\gamma$ lies on $\delta$, so its local time at 0 provides an additive functional that measures the size of the right boundary. Then we can translate the excursion decomposition of the Bessel process $X$ away from 0 into a decomposition of the $\text{SLE}_r(\kappa, \kappa-4)$ process in excursions away from $\delta$. We also recover the dimension of frontier times, proved to be equal to $1/2 + 2/\kappa$ in [3]. It is also possible to study a two-sided analogue of this situation, corresponding to the decomposition for $\text{SLE}(\kappa, \kappa - 4, \kappa - 4)$ starting from $(0, 0^-, 0^+)$ in excursions away from its cutpoints. This is closely related to the bead decomposition for Brownian (half-plane) excursions detailed in [24].

**Acknowledgments.** I wish to thank Wendelin Werner for his help and advice throughout the preparation of this paper. I also wish to thank Marc Yor for stimulating conversations Jean Bertoin and the Associate Editor for his careful rereading.

COURANT INSTITUTE
NEW YORK UNIVERSITY
251 MERCER STREET
NEW YORK, NEW YORK 10012
USA
E-MAIL: dubedat@cims.nyu.edu